\documentclass[a4paper,12pt,reqno]{amsart}
\usepackage{amssymb}
\usepackage{amsmath}
\usepackage{ifthen}
\usepackage[dvips]{graphicx}
\nonstopmode \numberwithin{equation}{section}
\usepackage{amssymb}
\usepackage{ifthen}
\usepackage{graphicx}
\usepackage{amsmath}
\usepackage[T1]{fontenc} %skandit
\usepackage{comment}

\nonstopmode \numberwithin{equation}{section}
\theoremstyle{plain}
\newtheorem{thm}{Theorem}
\newtheorem{cor}{Corollary}
\newtheorem{lem}[equation]{Lemma}
\newtheorem{prop}{Proposition}

\newtheorem{conj}{Conjecture}

\newcommand\numberthis{\addtocounter{equation}{1}\tag{\theequation}}
\theoremstyle{definition}
\newtheorem{defn}{Definition}[section]

\newtheorem{prob}{Problem}
\newtheorem{rem}{Remark}[section]

\newcounter{minutes}
\divide\time by 60
\newcounter{hours}
\multiply\time by 60
\addtocounter{minutes}{-\time}
\newcounter {own}
\def\theown {\thesection       .\arabic{own}}

\newenvironment{pf}[1][]{%
	\vskip 3mm
	\noindent
	\ifthenelse{\equal{#1}{}}%
	{{\slshape Proof. }}%
	{{\slshape #1.} }%
}%
{\qed\bigskip}

\newtheorem{Thm}{Theorem}

\newtheorem{Lem}{Lemma}

\newcommand{\ID}{{\mathbb D}}

\newcommand{\IC}{{\mathbb C}}

\renewcommand{\theequation}{\thesection.
\arabic{equation}}
\numberwithin{equation}{section}
\def\be{\begin{equation}}
\def\ee{\end{equation}}

\newcommand{\bee}{\begin{enumerate}}
	\newcommand{\eee}{\end{enumerate}}

\newcommand{\blem}{\begin{lem}}
	\newcommand{\elem}{\end{lem}}
\newcommand{\bthm}{\begin{thm}}
	\newcommand{\ethm}{\end{thm}}
\newcommand{\bcor}{\begin{cor}}
	\newcommand{\ecor}{\end{cor}}
\newcommand{\beg}{\begin{examp}}
	\newcommand{\eeg}{\end{examp}}
\newcommand{\begs}{\begin{examples}}
	\newcommand{\eegs}{\end{examples}}

\newcommand{\bdefn}{\begin{defn}}
	\newcommand{\edefn}{\end{defn}}

\newcommand{\bprob}{\begin{prob}}
	\newcommand{\eprob}{\end{prob}}
\newcommand{\bei}{\begin{itemize}}
	\newcommand{\eei}{\end{itemize}}

\newcommand{\bcon}{\begin{conj}}
	\newcommand{\econ}{\end{conj}}
\newcommand{\bcons}{\begin{conjs}}
	\newcommand{\econs}{\end{conjs}}
\newcommand{\bprop}{\begin{prop}}
	\newcommand{\eprop}{\end{prop}}
\newcommand{\br}{\begin{rem}}
	\newcommand{\er}{\end{rem}}
\newcommand{\brs}{\begin{rems}}
	\newcommand{\ers}{\end{rems}}
\newcommand{\bo}{\begin{obser}}
	\newcommand{\eo}{\end{obser}}
\newcommand{\bos}{\begin{obsers}}
	\newcommand{\eos}{\end{obsers}}
\newcommand{\bpf}{\begin{pf}}
	\newcommand{\epf}{\end{pf}}
\newcommand{\ba}{\begin{array}}
	\newcommand{\ea}{\end{array}}
\newcommand{\beq}{\begin{eqnarray}}
\newcommand{\beqq}{\begin{eqnarray*}}
\newcommand{\eeq}{\end{eqnarray}}
\newcommand{\eeqq}{\end{eqnarray*}}

\begin{document}

\title{The Landau-Bloch type theorems for certain class of holomorphic and pluriharmonic mappings in $\mathbb{C}^n$}

\author{Vasudevarao Allu}
\address{Vasudevarao Allu,
	School of Basic Science,
	Indian Institute of Technology Bhubaneswar,
	Bhubaneswar-752050, Odisha, India.}
\email{avrao@iitbbs.ac.in}

\author{Rohit Kumar}
\address{Rohit Kumar,
	School of Basic Science,
	Indian Institute of Technology Bhubaneswar,
	Bhubaneswar-752050, Odisha, India.}
\email{rohitk12798@gmail.com}

\subjclass[{AMS} Subject Classification:]{Primary 32A10, 31C10; Secondary 32A18, 30C62, 31B05.}

\def\thefootnote{}
\footnotetext{ {\tiny File:~\jobname.tex,
		printed: \number\year-\number\month-\number\day,
		\thehours.\ifnum\theminutes<10{0}\fi\theminutes }
} \makeatletter\def\thefootnote{\@arabic\c@footnote}\makeatother

\begin{abstract}
 In this paper, we first define two classes of holomorphic mappings defined on the unit ball $B^n$ of $n$ dimensional complex space $\IC^n$ and obtain the lower estimates for Bloch's constant for these classes. Also, we derive the Landau-Bloch type theorem for some subclasses of pluriharmonic mappings defined on the unit ball $B^n$.
% We find that if $f$ is a Bochner $K$-mapping defined on the unit ball $B^n$ into $\IC^n$ then $f(B^n)$ contains a schlicht ball of radius atleast $1/12K^{2n-1}.$ 
%
\end{abstract}

\maketitle
\pagestyle{myheadings}
\markboth{Vasudevarao Allu and  Rohit Kumar}{ The Landau-Bloch type theorems for certain class of holomorphic...}

\section{\textbf{Introduction}}
 Let $\mathbb{D}=\{z\in\mathbb{C}:|z|<1\}$ be the unit disc in the complex plane $\mathbb{C}$. In the case of one complex variable, the following theorem of Bloch is well known (see \cite{Bloch-1925}).

\begin{Thm}\cite{Bloch-1925}
Let $f$ be a holomorphic function on the $\overline{\mathbb{D}}=\{z:|z|\leq 1\}$ and satisfying $|f'(0)|=1$. Then there exists a positive constant $b$ such that $f(\mathbb{D})$ contains a schlicht disc of radius $b$.
\end{Thm}

 By a schlicht disc, we mean a disc which is the univalent image of some region in the unit disc $\mathbb{D}$. Let $\beta_f$ denote the least upper bound of the radii of all schlicht discs that $f$ carries and $\mathcal{F}$ denote the set of all holomorphic functions defined on $\overline{\mathbb{D}}=\{z:|z|\leq 1\} $ satisfying $|f'(0)|=1$, then the Bloch constant is the number defined by $$\beta(\mathcal{F})=\inf\{{\beta_f: f \in \mathcal{F}}\}.$$

In 1929, Landau \cite{Landau1} proved that if we replace the holomorphicity condition on $|z|\leq 1$ to $|z|< 1$, then also the corresponding constant is same. If one consider the function $f(z)=z$, then clearly $\beta(\mathcal{F})\leq1$. However, better estimates than these are known. The exact value of $\beta(\mathcal{F})$ is still unknown. 

In 1938, Alhforse and Grunsky \cite{Ahlfors-Grunsky-1937} proved that 
$$
0.4330 \approx \frac{\sqrt{3}}{4} \leq \beta(\mathcal{F}) \leq \frac{1}{\sqrt{1+\sqrt{3}}}\frac{\Gamma(1/3)\Gamma(11/12)}{\Gamma(1/4)}\approx 0.4719.
$$
 It is conjectured that this upper bound is the precise value of the Bloch constant and this conjecture is known as the Ahlfors-Grunsky conjecture which is still open. In 1990, Bonk \cite{Bonk-1990} obtained a slight improvment of the lower bound as $\beta(\mathcal{F})>\sqrt{3}/4+10^{-14}$. In 1996, Chen and Gauthier \cite{Bconstant} improved the result of Bonk and  proved that $ \beta(\mathcal{F})> \sqrt{3}/4+2\times10^{-4}$.
\vspace{2mm}

Let $$\mathbb{C}^n=\{z=(z_1,z_2,...,z_n):z_1,z_2,...,z_n\in\mathbb{C}\}$$
be the complex space of dimension $n$. For any point $z=\left(z_1, z_2, \cdots, z_n\right) \in \mathbb{C}^n$, 
$$
|z|=\left(\left|z_1\right|^2+\left|z_2\right|^2+\cdots+\left|z_n\right|^2\right)^{1 / 2}.
$$
We denote a ball in $\mathbb{C}^n$ with center at $z_0$ and radius $r$ by
$$
B^n\left(z_0, r\right)=\left\{z \in \mathbb{C}^n:\left|z-z_0\right|<r\right\}
$$
and the unit ball in $\mathbb{C}^n$ by
$$
B^n=\{z \in \mathbb{C}^n:|z|<1\}.
$$
For a mapping $f=\left(f_1, f_2, \cdots, f_n\right)$ of a domain in $\mathbb{C}^n$ into $\mathbb{C}^n$, we denote by $\partial f / \partial z_k$ the column vector formed by $\partial f_1 / \partial z_k, \partial f_2 / \partial z_k, \cdots, \partial f_n / \partial z_k$, and we denote by
$$
f^{\prime}=\left(\frac{\partial f}{\partial z_1}, \frac{\partial f}{\partial z_2}, \cdots, \frac{\partial f}{\partial z_n}\right),
$$
the matrix formed by these column vectors.
For an $n\times n$ matrix $A$, we have the matrix norm
$$
\|A\|=\left(\sum_{i, j}\left|a_{i j}\right|^2\right)^{1 / 2}
$$
and the operator norm 
$$|A|=\sup_{z\neq 0}\frac{|Az|}{|z|}.$$

In the case of several complex variables, the classical theorem of Bloch for holomophic mappings in the disc fails to extend to general holomorphic mappings in the ball of $\mathbb{C}^n$. In 1967, Wu \cite{Wu-1967}  pointed out that the Bloch theorem fails unless we have some restrictive assumption on holomorphic mappings. For example, for positive integer $n$, setting $ f_n(z_1, z_2) = (nz_1,z_2/n)$, we see that $f_n$ is an holomorphic function on the unit ball of $\mathbb{C}^2$ for each $n$ and $|detJ_{f_n}(0)|= 1$. But the image of $f$ contains the schlicht ball of radius at most $1/n$. So, the infimum of the radii of the schlicht balls
anywhere when all members of $\{f_n\}$ are taken into account is zero. In order that there is a positive Bloch constant, it is necessary to restrict the class of holomorphic mappings in higher dimension. Earlier investigations for some subclass of holomorphic mappings were done by Bochner \cite{Bochner-1946}, Hahn \cite{Hahn-1973}, Harris \cite{Harris-1977}, Sakaguchi \cite{Sakaguchi-1956}, Takahashi \cite{Takahashi-1951}, and Wu \cite{Wu-1967}. If we drop the holomorphicity, there is no Bloch theorem for quasiregular mappings of the ball, but Eremenko \cite{Eremenko-2000} has proved that a Bloch theorem for entire quasiregular mappings.
\vspace{2mm}

In 1946, Bochner \cite{Bochner-1946} defined a subclass of the class of holomorphic functions from $B^n\left(\subseteq \mathbb{C}^n\right)$ to $\mathbb{C}^n$. For a constant $K \geq 1$, a holomorphic mapping $f$ from $B^n$ to $\mathbb{C}^n$ is said to be Bochner $K$-mapping in $B^n$ if $f$ satisfies the following differential inequality
$$
\left\|f'\right\| \leq K\left|\operatorname{det} f'\right|^{1/n}
$$
at each $z \in B^n$. Bochner \cite{Bochner-1946} has proved  that, for each $K \geq 1$ and $n \geq 2$, there is a constant $\beta>0$ such that, for each normalized Bochner $K$-mapping $f$ in the unit ball, $\beta_f \geq \beta$ holds. However, Bochner has not given an estimate for this Bloch constant.
\vspace{2mm}

In 1951, Takahashi \cite{Takahashi-1951} defined a new class of normalized holomorphic mappings $f$ satisfying the weaker condition
\begin{equation}\label{Rohi-Vasu-P2-equation-000}
\max_{|z|<r}\|f'(z)\|\leq K\max_{|z|\leq r}|\operatorname{det} f'(z)|^{1/n}, \quad \text{for each}\ 0\leq r<1.
\end{equation}
The holomorphic functions which satisfies (\ref{Rohi-Vasu-P2-equation-000}) are called Takahashi $K$-mappings. For such normalized Takahashi $K$-mappings, Takahashi \cite{Takahashi-1951} has proved that
$$\beta\geq \frac{(n-1)^{n-2}}{12K^{2n-1}}.$$
\vspace{2mm}
Later in 1956, Sakaguchi \cite{Sakaguchi-1956} improved Takahashi's estimate to
$$\beta \geq \frac{(n-1)^{n-2}}{8K^{2n-1}}.$$ 
\vspace{2mm}

We now define a new class of holomorphic mappings from $B^n$ into $\mathbb{C}^n$ which contains the Bochner $K$-mappings, we call such functions as Bochner $\left(K, K'\right)$- mappings.
\begin{defn}\label{Rohi-Vasu-P2-Definition-001}
For constants $K\geq1, K'\geq0$, a holomorphic mapping $f$ from $B^n$ into $\mathbb{C}^n$ is said to be Bochner $\left(K, K'\right)$- mapping in $B^n$ if $f$ satisfies the following differential inequality
$$
\left\|f'\right\|^2 \leq K^2\left|\operatorname{det} f'\right|^{2/n}+K' 
$$
at each $z \in B^n$.
\end{defn}

We remark that unit ball $B^n$ in the definition \ref{Rohi-Vasu-P2-Definition-001} can be replaced by any general domain in $\mathbb{C}^n$. In particular, if $K'=0$, then Bochner $(K,K')$-mapping reduces to Bochner $K$-mapping. We see that every Bochner $K$-mapping is a Bochner $(K,K')$-mapping for $K'=0$, but the converse need not be true. This can be seen from the following example: Let $$f(z_1,z_2)=\left(z_1+z_2, \left(z_1-\frac{1}{2}\right)^2+\left(z_2-\frac{1}{3}\right)^2\right)$$ in $B^2$, which is clearly holomorphic mapping on $B^2.$ One can easily see that the mapping $f(z_1,z_2)$ is not Bochner $K$-mapping for any $K\geq 1$ but it is a Bochner $(1,20)$- mapping.
\vspace{2mm}

 The motivation for defining this type of mappings came from the paper of Nirenberg (see \cite{Nirenberg-1953}), where Nirenburg has defined this type of mappings in the plane. For more details on this type of mappings, we refer to  \cite{Allu-Kumar-2024,Chen-Li-Sahoo-Wang-2017,Chen-Ponnusamy-2020,Chen-Ponnusamy-Wang-2021,Finn-Serrin-1958,Nirenberg-1953}.
 \vspace{2mm}

Let $\lambda_f^2(z)$ and $\Lambda_f^2(z)$ denote the smallest and the largest eigenvalues of the Hermitian matrix $A^*A$, where $A=f'(z)$ and $A^*$ is the conjugate of $A$. A holomorphic mapping $f$ from the unit ball $B^n$ of $\IC^n$ into $\IC^n$ is $K$-quasiregular if
 $$\Lambda_f(z)\leq K \lambda_f(z)$$ 
 at every point $z \in B^n$. A mapping is said to be quasiregular if it is $K$-quasiregular for some $K\geq 1$. From \cite{Marden-Rickman-1974} we know that for $n>1$, a quasiregular holomorphic mappings are locally biholomorphic . In fact, Poletsky \cite{Poletsky-1985} also proved that quasiregular holomorphic mappings(in any bouded domain) are rather rigid.  In 2000, Chen and Gauthier \cite{Chen-Gauthier-2000} proved that if $f$ is a $K$- quasiregular holomorphic mapping with the normalization $\operatorname{det} f'(0)=1$, then the image $f(B^n)$ contains a schlicht ball of radius at least $1/12K^{1-1/n}$.
 \vspace{2mm}

Now, we define a new class of holomorphic mappings from $B^n$ into $\mathbb{C}^n$ which contains the $K$-quasiregular mappings, we call such mappings $\left(K, K'\right)$- quasiregular mappings.
\begin{defn}\label{Rohi-Vasu-P2-Definition-002}
For constants $K\geq1, K'\geq0$, a holomorphic mapping $f$ from $B^n$ into $\mathbb{C}^n$ is said to be $\left(K, K'\right)$-quasiregular mapping in $B^n$ if 
$$
\Lambda_f(z)\leq K\lambda_f(z)+K' 
$$
at each $z \in \mathbb{B}^n$.
\end{defn}
A continuous complex-valued function $\phi$ defined on a domain $\Omega \subset \IC^n$ is called a pluriharmonic mapping if, for each fixed $z'\in\Omega$ and $\theta \in \partial B^n$, the function $\phi\left(z'+\theta \zeta\right)$ is harmonic in $\{\zeta : |\zeta| < d_\Omega(z)\}$, where $d_\Omega(z)$ denotes the distance from $z$ to the boundary $\partial\Omega$ of $\Omega.$ A mapping $f$ of $\Omega$ into $\IC^n$ is called a pluriharmonic mapping if every component of $f$ is pluriharmonic. A mapping $f$ of $B^n$ into $C^n$ is pluriharmonic if, and only if, $f$ has a representation $f=g+\bar{h}$, where $g$ and $h$ are holomorphic mappings (see \cite{Rudin-1980}). 
\vspace{2mm}

For a continuously differentiable mapping $w=f(z)=\left(f_1(z), \ldots f_m(z)\right): B^n \rightarrow \mathbb{C}^m, z=\left(z_1 \ldots \ldots z_n\right)$, by $f_z$ and $f_{\bar{z}}$ denote the matrices $\left(\partial f_j / \partial z_k\right)_{m \times n}$ and $\left(\partial f_j / \partial \bar{z}_k\right)_{m \times n}$, respectively.
Denote the maximum dilation $\widetilde{\Lambda}_f$ and minimum dilation $\widetilde{\lambda}_f$ by
$$
\widetilde{\Lambda}_f(z)=\max _{\theta \in \partial B^n}\left|f_z(z) \theta+f_{\bar{z}}(z) \bar{\theta}\right|\quad\text{and }\widetilde{\lambda}_f(z)=\min _{\theta \in \partial B^n}\left|f_z(z) \theta+f_{\bar{z}}(z) \bar{\theta}\right|,
$$
respectively, where $\theta$ is regarded as a column vector. Also, when $n=1$, it is easy to see the corresponding definitions for planar harmonic mappings.
\vspace{2mm}

Wang {\it et al.} \cite{Wang-Yang-Liu-2022} generalized the notion of $K$-quasiregular mapping to pluriharmonic mappings and established a lower bound estimate of Bloch constant for such mappings.
\vspace{2mm}

A pluriharmonic mapping $f$ of $B^n$ into $\IC^n$ is said to be $K$-quasiregular pluriharmonic if
 $$
       \widetilde{\Lambda}_f(z)\leq K\widetilde{\lambda}_f(z)^{1/n}\quad \text{for} \ z\in B^n. 
 $$
Now, we define a new class of pluriharmonic mappings of $B^n$ into $\mathbb{C}^n$ which contains $K$-quasiregular pluriharmonic mappings, we call such mappings $\left(K, K'\right)$- quasiregular pluriharmonic mappings.
\begin{defn}\label{Rohi-Vasu-P2-Definition-003}
Let $f: B^n \to \IC^n$ be a pluriharmonic mapping and $K\geq 1, K'\geq 0$. We say that $f$ is a $(K,K')$-quasiregular pluriharmonic mapping if $$\widetilde{\Lambda}_f(z)\leq K\widetilde{\lambda}_f(z)^{1/n}+K'$$ for $z\in B^n$.

In particular, if $K'=0$, then $(K,K')$-quasiregular pluriharmonic mapping reduces to $K$-quasiregular pluriharmonic mapping. 
\end{defn}
\vspace{2mm}
In 2011, Chen and Gauthier \cite{Chen-Gauthier-2011} established the Landau theorems and Bloch theorems for pluriharmonic mappings $f: B^n \to \IC^n.$ S. Chen, S. Ponnusamy and X. Wang  have studied Landau-Bloch constants for some specified spaces such as, pluriharmonic Bergman space, $\alpha$-Bloch space and hyperbolic-harmonic Bloch space (see \cite{Chen-Ponnusamy-Wang-2012-a,Chen-Ponnusamy-Wang-2012-b,Chen-Ponnusamy-Wang-2012-c,Chen-Ponnusamy-Wang-2015}). In 2020, Xu and Liu  \cite{Xu-Liu-2020} obtained a new version of the Bloch theorem for pluriharmonic $\nu$ -Bloch-type mappings. In 2022, Liu and Ponnusamy \cite{Liu-Ponnusamy-2022} obtained three Bloch-type theorems of pluriharmonic mappings in $B^n$, which improve the corresponding results of Chen and Gauthier \cite{Chen-Gauthier-2011}.

\section{ The Bloch theorem for Bochner $(K,K')$-mappings}
The following result of Takahashi \cite{Takahashi-1951} is useful to prove our main result in this section.

\begin{Thm}\cite{Takahashi-1951}\label{Takahashi}
Let $f(z)=\left(f_1(z),f_2(z),...,f_n(z)\right)$ be an analytic transformation defined by $n$ functions $f_i(z)$ of $n$ complex variables $z$, each analytic in a domain $D$ in $n$ dimensional complex space of the variables $z$, and let its Jacobian $J_f(z)$ does not vanish at a point $a$ in $D$. Let $\partial f_i(a) / \partial z_k \equiv \alpha_{i k} ; i, k=1, \cdots, n ; A=\left(\alpha_{i k}\right)$, $\operatorname{det} A=J_f(a) \neq 0$, it follows that the characteristic values $\lambda_1, \cdots, \lambda_n$ of $A^*A$ are real and positive, so that $\min \{\lambda_1, \cdots, \lambda_n\}=\lambda>0$. Next, let $\rho_0$ be the upper limit of $\rho$ such that the inequality
$$
\sum_{i, k=1}^n\left|\frac{\partial f_i}{\partial z_k}(z)-\frac{\partial f_i}{\partial z_k}(a)\right|^2 \leq \lambda
$$
is satisfied for $\sum_{k=1}^{n} |z_k-a_k|^2\leq \rho$. Then $f$ is univalent on a ball with center $a$ and radius $\rho_0$, i.e., $B^n(a,\rho_0)$. Further, $f(B^n(a,\rho_0))$ contains a ball with center $f(a)$ and radius $2^{-1} \lambda^{\frac{1}{2}} \rho_0$. Moreover, the value $2^{-1} \lambda^{\frac{1}{2}} \rho_0$ can not be replaced by larger one for certain analytic transformation and for some $a$.
\end{Thm} 
\vspace{1.5mm}

The following lemma gives an estimates the lower bound of the smallest singular values of a non-singular matrix.

\begin{Lem}\label{lemma1}\cite{Yu-Gu-1997}
 If $A$ is a non-singular $n \times n$ matrix. Then $A^*A$ is a positive definite hermitian matrix, the characteristic values of $\lambda_1, \lambda_2, \ldots, \lambda_n$ of $A^* A$ are therefore real and positive, so that $\lambda=\min \left\{\lambda_2, \ldots, \lambda_n\right\}>0$. Then the following inequality is satisfied
$$
\lambda>(n-1)^{n-1}|\operatorname{det} A|^2\|A\|^{-2(n-1)},
$$
$\|A\|$ is the Euclidean norm of the matrix A.
\end{Lem}

\vspace{2mm}
In this section, we show that for each $K \geq 1, \ K' \geq 0$ and $n \geq 2$, there is a constant $\beta>0$ such that, for each normalized Bochner $\left(K, K^{\prime}\right)$-mapping $f$ in the ball $\beta_f \geq \beta$.

\begin{thm}\label{Rohi-Vasu-P2-Theorem-001}
  Let $n \geq 2$ be any integer and $K \geq 1, \ K^{\prime} \geq 0$ be any  constants. Let $ f(z)$ be Bochner $\left(K, K'\right)$- mapping from $\overline{B}^n$ to $\mathbb{C}^n$ with $|\operatorname{det} f'(0)|=1$, then $f$ maps some subdomain of unit ball univalently onto a ball $B^n$ of positive radius 
$$R(n,K,K')=\frac{1}{4(K^2+K')^{n-1}}\left(\sqrt{4K^2+K'}+\sqrt{K^2+K'}\right)^{-1}.$$
In other words, we say that Bloch Theorem holds for Bochner $(K, K')$- mappings.
\end{thm}
\begin{pf} [{\bf Proof}]
Let $f(z)=(f_1(z), f_2(z),...,f_n(z))$ be Bochner $(K,K')$- mapping from $\overline{B}^n$ to $\mathbb{C}^n$ with $|det f'(0)|=1$. We introduce the functions
$$
M(r)=  \max _{|z| \leq r}\left|\operatorname{det} f^{\prime}(z)\right|^{1/n}, 
$$ and  $$  \phi(r)=r M(1-r),
$$
for $0 \leq r \leq 1 $.\\
We see that
$$
\phi(0)=0\ \text{and}\ \phi(1)=M(0)=\left|\operatorname{det} f^{\prime}(0)\right|^{1/n}=1.
$$
Then there exists a number $r_0, \ 0<r_0 \leq 1$ such that
$
\phi\left(r_0\right)=1
$
and $\phi(r)<1$ for $0 \leq r<r_0.$
This implies
$$
M\left(1-r_0\right)=\frac{\phi\left(r_0\right)}{r_0}=\frac{1}{r_0},
$$
and for any $0<r<r_0$,
\begin{equation}\label{Rohi-Vasu-P2-equation-001}
M(1-r)=\frac{\phi(r)}{r}<\frac{1}{r}. 
\end{equation} 
Let $\alpha=(\alpha_1,\alpha_2,...,\alpha_n)$ be any point inside the closed ball of radius $1-r_0$ such that
 $$\left|\operatorname{det} f^{\prime}(\alpha)\right|^{1 / n}=M\left(1-r_0\right)=\frac{1}{r_0}.$$
Define $ F:\overline{B}^n \rightarrow \mathbb{C}^n $ by
$$
 F(\zeta)=F\left(\zeta_1, \zeta_2, \ldots, \zeta_n\right)=(F_1(\zeta),F_2(\zeta),\cdots,F_n(\zeta))=2\left(f\left(\alpha+\frac{r_0}{2} \zeta\right)-f(\alpha)\right) 
$$
for $|\zeta| \leq 1.$
Then clearly $F(0)=0$.

\vspace{2mm}
It is easy to see that  for $|\zeta| \leq 1$,
$$
\left|\alpha+\frac{r_0}{2} \zeta\right| \leq|\alpha|+\frac{r_0}{2}|\zeta| \leq 1-r_0+\frac{r_0}{2}=1-\frac{r_0}{2}<1.
$$
Therefore, $F$ is well-defined.\\
A simple computaion shows that
$$
\begin{aligned}
\frac{\partial F_i}{\partial \zeta_j}(\zeta) & =2\frac{\partial f_i}{\partial \zeta_j}\left(\alpha+\frac{r_0}{2} \zeta\right) \\ \vspace{2mm}
& =r_0\frac{\partial f_i}{\partial z_j}\left(\alpha+\frac{r_0}{2} \zeta\right)
\end{aligned}
$$
where $z_j=\alpha_j+\frac{r_0}{2}\zeta_j, \quad j=1,2,...,n.$\\
Therefore,$$ \operatorname{det} F'(\zeta)={r_0}^n \operatorname{det} f'\left(\alpha+\frac{r_0}{2} \zeta\right), 
$$
which implies that
$$
\left|\operatorname{det} F^{\prime}(0)\right|={r_0}^n\left|\operatorname{det} f^{\prime}(\alpha)\right|={r_0}^n\left(M\left(1-r_0\right)\right)^n={r_0}^n\left(\frac{1}{r_0}\right)^n 
=1.
$$
A simple computation shows that
$$
\begin{aligned}
\left\|F'(\zeta)\right\|^2 & =\sum_{i, j}\left|\frac{\partial F_i}{\partial \zeta_j}(\zeta)\right|^2=\sum_{i, j}{r_0}^2\left|\frac{\partial f_i}{\partial z_j}\left(\alpha+\frac{r_0}{2} \zeta\right)\right|^2 \\
& ={r_0}^2 \sum_{i, j}\left|\frac{\partial f_i}{\partial z_j}\left(\alpha+\frac{r_0}{2} \zeta\right)\right|^2 \\
& ={r_0}^2\left\|f'\left(\alpha+\frac{r_0}{2} \zeta\right)\right\|^2 \\
& \leq {r_0}^2\left[K^2\left|\operatorname{det} f'\left(\alpha+\frac{r_0}{2} \zeta\right)\right|^{2 / n}+K'\right].
\end{aligned}
$$
Since
$\left|\alpha+\frac{r_0}{2} \zeta\right| \leq 1-r_0/2$. By the definition of the function $M$, we obtain
$$
\left|\operatorname{det} f'\left(\alpha+\frac{r_0}{2} \zeta\right)\right|^{2 / n} \leq\left(M\left(1-\frac{r_0}{2}\right)\right)^2.
$$
From (\ref{Rohi-Vasu-P2-equation-001}), it follows that
\begin{align*}\label{Rohi-Vasu-P2-equation-002}
\left\|F'(\zeta)\right\|^2 &\leq {r_0}^2\left(K^2\left(M\left(1-\frac{r_0}{2}\right)\right)^2+K'\right) \\
&\leq{r_0}^2\left(K^2  \frac{4}{{r_0}^2}+K'\right)\\
& =4K^2+{r_0}^2 K' \\
& \leq 4K^2+K'.\numberthis
\end{align*}

Further,

\begin{align*}\label{Rohi-Vasu-P2-equation-003}
\left\|F^{\prime}(0)\right\|^2=&{r_0}^2 \sum_{i, j}\left|\frac{\partial f_i}{\partial z_j}(\alpha)\right|^2 \\ &={r_0}^2\left\|f^{\prime}(\alpha)\right\|^2 \\
& \leq {r_0}^2\left(K^2\left|\operatorname{det} f'(\alpha)\right|^{2/n}+K'\right) \\
& \leq {r_0}^2\left(K^2 \cdot \frac{1}{r_0^2}+K'\right) \\
&\leq K^2+K'.\numberthis
\end{align*}
\vspace{2mm}

From (\ref{Rohi-Vasu-P2-equation-002}) and (\ref{Rohi-Vasu-P2-equation-003}), it is easy to see that
\begin{equation}\label{Rohi-Vasu-P2-equation-004}
\left\|F'(\zeta)-F'(0)\right\| \leq\left\|F'(\zeta)\right\|+\left\|F'(0)\right\| \leq \sqrt{4K^2+K'}+\sqrt{K^2+K'}
\end{equation}

for $\left|\zeta\right|\leq 1$.
\vspace{2mm}

By applying the Schwarz lemma for functions of several complex variables (see \cite{Bochner-Martin-1948}) to (\ref{Rohi-Vasu-P2-equation-004}), we obtain
\begin{equation}\label{Rohi-Vasu-P2-equation-004a}
\left\|F'(\zeta)-F'(0)\right\| \leq \left(\sqrt{4K^2+K'}+\sqrt{K^2+K'}\right)|\zeta|
\end{equation}
for $|\zeta|\leq 1.$ Since $\left|\operatorname{det} F^{\prime}(0)\right|=1 \neq 0$, the characteristic values $\lambda_1, \lambda_2, \ldots, \lambda_n$ of $\left(F^{\prime}(0)\right)^* F^{\prime}(0)$ are real and positive, so that $\min \left\{\lambda_1, \lambda_2, \ldots, \lambda_n\right\}=\lambda>0$, and the inequality (\ref{Rohi-Vasu-P2-equation-004a}) takes the following form
\begin{equation}\label{Rohi-Vasu-P2-equation-005}
\left\|F'(\zeta)-F'(0)\right\|^2 \leq \lambda \quad \text { for }\ |\zeta| \leq \frac{\lambda^{1 / 2}}{\sqrt{4K^2+K'}+\sqrt{K^2+K'}}.
\end{equation}
Since the product of eigenvalues is equal to determinant of the matrix, we have
$$
 \lambda^n \leq \lambda_1 \lambda_2 \cdots \lambda_n=|\operatorname{det} F'(0)|^2. $$
 This implies 
 $$
 \lambda^{1 / 2} \leq\left|\operatorname{det} F^{\prime}(0)\right|^{1 / n}=1
$$
and hence
\begin{equation}\label{Rohi-Vasu-P2-equation-006}
  \frac{{\lambda}^{1/2}}{\sqrt{4K^2+K'}+\sqrt{K^2+K'}}\leq \frac{1}{\sqrt{4K^2+K'}+\sqrt{K^2+K'}}<1.
\end{equation}

In view of (\ref{Rohi-Vasu-P2-equation-006}) together with inequality (\ref{Rohi-Vasu-P2-equation-005}) and applying Theorem \ref{Takahashi}, we see that $F(\zeta)$ maps some subdomain of $|\zeta| \leq 1$ univalently  onto a ball with center $0$ and radius
$$
\frac{\lambda^{1 / 2}}{2} \frac{{\lambda}^{1/2}}{\sqrt{4K^2+K'}+\sqrt{K^2+K'}}=\frac{ \lambda}{2(\sqrt{4K^2+K'}+\sqrt{K^2+K'})}.
$$
Using Lemma \ref{lemma1}, for $n\geq2$ we obtain the following
$$
\lambda >(n-1)^{n-1} \left|\operatorname{det} F^{\prime}(0)\right|^2\left\|F^{\prime}(0)\right\|^{-2(n-1)}\geq   \left|\operatorname{det} F^{\prime}(0)\right|^2\left\|F^{\prime}(0)\right\|^{-2(n-1)}.
$$
By (\ref{Rohi-Vasu-P2-equation-003}), we have 
$$\lambda>(K^2+K')^{-(n-1)}.
$$
Therefore, $F(\zeta)$ maps some sub-domain of the ball $|\zeta|\leq1$ univalently onto a ball with center $0$ and  of radius atleast
$$
\frac{1}{2}\frac{(K^2+K')^{-(n-1)}}{\sqrt{4K^2+K'}+\sqrt{K^2+K'}} .
$$
Hence by the definition of $F(\zeta)$, $f(z)$ maps some sub-domain of the unit ball univalently onto a ball of radius atleast
$$R(n,K,K')=\frac{1}{4(K^2+K')^{n-1}}\left(\sqrt{4K^2+K'}+\sqrt{K^2+K'}\right)^{-1}.$$
This completes the proof.
\end{pf}\\
Using Theorem \ref{Rohi-Vasu-P2-Theorem-001}, we obtain the following corollary.
\begin{cor}\label{Rohi-Vasu-P2-Corollary-001}
 Let $n \geq 2$ be any integer and $K \geq 1$ be any positive constant. Let
$ f(z)$ be Bochner $K$- mapping from $\overline{B}^n$ to $\mathbb{C}^n$ with $|\operatorname{det} f'(0)|=1$, then $f$ maps some subdomain of unit ball $B^n$ univalently onto a ball of positive radius $$R(n,K)=\frac{1}{12K^{2n-1}}.$$
\end{cor}
\begin{pf} [{\bf Proof}]
Let $f(z)$ be a Bochner $K$-mapping. Then we can easily see that $f(z)$ is a Bochner $(K,0)$- mapping. Therefore, the proof follows by substituting $K'=0$ in the Theorem \ref{Rohi-Vasu-P2-Theorem-001}.
\end{pf}
\vspace{2mm}

\section{The Bloch theorem for $(K,K')$-quasiregular mappings}

Following lemma is useful in proving our main result in this section.
\begin{Lem}\cite{Hahn-1973}\label{lemma2}
Let $w=f(z)$ be a holomorphic mapping defined in a neighborhood of a point $t \in \IC^n$ into $\IC^n$ with $J_f(t) \neq 0$. Suppose that $\lambda_f \equiv \lambda_f(t)$ is the positive square root of the smallest characteristic value of the matrix $A^* A$ at $t$, where $A \equiv(d f / d z)$. Then the following hold:
\begin{enumerate}

\item The mapping $w=f(z)$ is univalent in any open convex subset $K, t \in K$, of the set
$$
\Omega_f=\left\{z:|A(z)-A(t)|<\lambda_f\right\},
$$
where $|A|=\sup _{|x|=1}|A x|$ and $|x|$ denotes the euclidean norm of the n-vector $x$.

\item If $r_0$ is the radius of the largest ball contained in $\Omega_f$ centered at $t$, then $f\left[B\left(t, r_0\right)\right]$ contains the ball of radius $r_0 \lambda_f / 2$ centered at $f(t)$, where $B\left(t, r_0\right)$ $=\left[z:|z-t|<r_0\right]$.
\end{enumerate}
\end{Lem}

In this section we show that for each $K \geq 1, \ K' \geq 0$ and $n \geq 2$, there is a constant $\beta>0$ such that, for each normalized $\left(K, K^{\prime}\right)$-quasiregular mapping $f$ in the ball $\beta_f \geq \beta$.

\begin{thm}\label{Rohi-Vasu-P2-Theorem-002}
Let $f:B^n\to \IC^n$ be a $(K,K')$-quasiregular mapping of the unit ball $B^n$ into $\IC^n$ with $\lambda_f(0)\geq \alpha>0$. Then
$$\beta_f\geq \frac{\alpha^2}{4(2K\alpha+K'+\alpha)}.$$ 
\end{thm}
\begin{pf}[{\bf Proof}]
Let $f:B^n\to \IC^n$ be a $(K,K')$-quasiregular mapping of the unit ball $B^n$ into $\IC^n$ with $\lambda_f(0)\geq \alpha>0$. Without loss of generality, we assume that $f$ is holomorphic on $\overline{B}^n$. We introdce the following functions
$$N(r)=\max_{|z|\leq r}\lambda_f(z)$$
and
$$\psi(r)=rN(1-r)$$
for $0\leq r\leq1.$
We see that
 $$ \psi(0)=0 \text{ and  }\psi(1)=N(0)=\lambda_f(0)=\alpha>0.$$
Then there exists a number $r_0, \ 0<r_0\leq 1$ such that $\psi(r_0)=\alpha$ and $\psi(r)<\alpha$ for $0\leq r_0. $ This implies 
      $$N(1-r_0)=\frac{\psi(r_0)}{r_0}=\frac{\alpha}{r_0},$$
and for any $0<r<r_0$,
\begin{equation}\label{Rohi-Vasu-P2-equation-007}
 N(1-r)=\frac{\psi(r)}{r}<\frac{\alpha}{r}.
\end{equation}
Let $w_0$ be any point inside the closed ball of radius $1-r_0$ such that
$$\lambda_f(w_0)=N(1-r_0)=\frac{\alpha}{r_0}.$$
Define $G:B^n\to \IC^n$ by
\begin{equation}\label{Rohi-Vasu-P2-equation-007-a}
G(\zeta)=\frac{r_0}{2}A^{-1}\left(f(w_0+\frac{r_0}{2}\zeta)-f(w_0)\right)
\end{equation}
for $|\zeta|\leq 1$ and $A=f'(w_0).$ It is easy to see that $$|w_0+\frac{r_0}{2}\zeta|\leq |w_0|+\frac{r_0}{2}|\zeta|\leq 1-r_0+\frac{r_0}{2}=1-\frac{r_0}{2}<1.$$
Also, $[f'(w_0)]^{-1}$ exists because $|\operatorname{det} f'(w_0)|\geq \lambda_f^n(w_0)>0.$ Therefore, $G$ is well defined.
\vspace{1mm}
\begin{flushleft}
We observe that 
\end{flushleft}
$$\frac{dG(\zeta)}{d\zeta}=\frac{r_0}{2}A^{-1}\frac{df}{d\zeta}\left(w_0+\frac{r_0}{2}\zeta\right), $$
which implies
\begin{equation}\label{Rohi-Vasu-P2-equation-007-b}
G'(\zeta)=A^{-1}f'(w_0+\frac{r_0}{2}\zeta).
\end{equation}
By using Cauchy–Schwarz inequality, we obtain
$$\Lambda_G(\zeta)=|G'(\zeta)|\leq |A^{-1}||f'(w_0+\frac{r_0}{2}\zeta)|.$$
Since $|A^{-1}|=\frac{1}{\lambda_f(w_0)}$, we have
$$\Lambda_G(\zeta)\leq \frac{\Lambda_f(w_0+\frac{r_0}{2}\zeta)}{\lambda_f(w_0)}.$$
Since $f$ is $(K,K')$-quasiregular mapping, we have 
$$\Lambda_G(\zeta)\leq \frac{K\lambda_f(w_0+\frac{r_0}{2}\zeta)+K'}{\lambda_f(w_0)}.$$
Since $|w_0+\frac{r_0}{2}\zeta|\leq 1-r_0/2.$ By the definition of the function $N$, we obtain
$$\lambda_f(w_0+\frac{r_0}{2}\zeta)\leq N\left(1-\frac{r_0}{2}\right).$$
Using (\ref{Rohi-Vasu-P2-equation-007}), we obtain
\begin{align*}\label{Rohi-Vasu-P2-equation-008}
\Lambda_G(\zeta)&\leq \frac{KN\left(1-\frac{r_0}{2}\right)+K'}{\lambda_f(w_0)}\\
&\leq \frac{2K\alpha/r_0+K'}{\alpha/r_0}\leq \frac{1}{\alpha}(2K\alpha+K'). \numberthis
\end{align*}
From (\ref{Rohi-Vasu-P2-equation-007-b}), we  have $G'(0)=I_n$ and hence, $\Lambda_G(0)=\lambda_G(0)=1.$ Thus,
\begin{equation}\label{Rohi-Vasu-P2-equation-009}
|G'(\zeta)-G'(0)|\leq |G'(\zeta)|+|G'(0)|\leq \frac{1}{\alpha}(2K\alpha+K')+1=\frac{1}{\alpha}\left(2K\alpha+K'+\alpha\right)
\end{equation}
for $|\zeta|\leq 1.$
By the Schwarz lemma (see \cite{Bochner-Martin-1948}), we obtain 
\begin{equation}\label{Rohi-Vasu-P2-equation-010}
|G'(\zeta)-G'(0)|\leq \frac{1}{\alpha}\left(2K\alpha+K'+\alpha\right)|\zeta|
\end{equation} 
for $|\zeta|\leq 1.$
Clearly (\ref{Rohi-Vasu-P2-equation-010}) shows that
\begin{equation*}
|G'(\zeta)-I_n|\leq 1 \quad \text{for} \quad |\zeta|\leq \frac{\alpha}{2K\alpha+K'+\alpha}.
\end{equation*}
By Lemma \ref{lemma2}, $w=G(\zeta)$ maps the ball $B^n\left(0,\alpha/(2K\alpha+K'+\alpha)\right)$ univalently onto a domain containing the ball $B^n\left(0,\alpha/(4K\alpha+2K'+2\alpha)\right).$ Hence by (\ref{Rohi-Vasu-P2-equation-007-a}), $w=f(z)$ maps the subdomain $B^n \left(w_0, r_0\alpha/(4K\alpha+2K'+2\alpha)\right)$ of $B^n$ univalently onto a ball center at $f(w_0)$ and radius 
$$\frac{\alpha^2}{4(2K\alpha+K'+\alpha)}.$$
This completes the proof.
\end{pf}

\begin{cor}
Let $f:B^n\to \IC^n$ be a $(K,K')$-quasiregular mapping of the unit ball $B^n$ into $\IC^n$ with $|\operatorname{det} f'(0)|=\alpha>0$. Then
$$\beta_f\geq \frac{\alpha^{2/n}}{4(2K\alpha^{1/n}+K'+\alpha^{1/n})}.$$ 
\end{cor}
\begin{pf}[{\bf Proof}]
Since $$\alpha=|\operatorname{det} f'(0)|\leq \lambda_f^n(0),$$
we have,
     $$\lambda_f(0)\geq\alpha^{1/n}.$$
Now by replacing $\alpha$ by $\alpha^{1/n}$ in Theorem \ref{Rohi-Vasu-P2-Theorem-002}, we obtain the desired result.
\end{pf}

\section{The Landau-Bloch type theorem for pluriharmonic mappings}

In 2011, Chen and Gauthier \cite{Chen-Gauthier-2011} proved the following Schwarz-Pick lemma for pluriharmonic mappings:
\begin{Lem}\label{lemma3}\cite{Chen-Gauthier-2011}
Let $f$ be a pluriharmonic mapping of $B^n$ into $B^m$. Then 
$$
\widetilde{\Lambda}_f(z)\leq \frac{4}{\pi}\frac{1}{1-|z|^2}\quad \text{for}\ z\in B^n.
$$
If $f(0)=0,$ then
$$
|f(z)|\leq \frac{4}{\pi}\arctan|z|\leq \frac{4}{\pi}|z|\quad\text{for} \ z\in B^n.
$$
\end{Lem}

The following Landau-type theorem for pluriharmonic mappings of $B^n$ into $\IC^n$ with bounded dilation has been proved by Wang {\it{et al.}} \cite{Wang-Yang-Liu-2022}.
\begin{Thm}\cite{Wang-Yang-Liu-2022}\label{Wang-2022}
Let $f$ be a pluriharmonic mapping of $B^n$ into $\IC^n$ such that $f(0)=0, \widetilde{\lambda}_f(0)=1$ and $\widetilde{\Lambda}_f(z)\leq \widetilde{\Lambda}$ for $z\in B^n$. Then $f$ is univalent on the ball $B^n(o,\rho)$ and the range $f(B^n(0,\rho))$ covers the ball $B^n(0,R)$, where
$$\rho=\frac{\pi}{4(\widetilde{\Lambda}_f(0)+\widetilde{\Lambda})}\quad\text{and}\ R=\frac{\pi}{8(\widetilde{\Lambda}_f(0)+\widetilde{\Lambda})}.$$
If, in addition, $\widetilde{\Lambda}_f(0)=1,$ then $f$ is univalent on the ball $B^n(0,\rho')$ and range $f(B^n(0,\rho'))$ covers the ball $B^n(0,R')$, where
$$\rho'=\frac{\pi}{4(1+\widetilde{\Lambda})}\quad\text{and}\ R'=\frac{\pi}{8(1+\widetilde{\Lambda})}.$$
\end{Thm}

For $(K,K')$-quasiregular pluriharmonic mappings with bounded dilation, we prove the following Landau-type theorem.
\begin{thm}\label{Rohi-Vasu-P2-Theorem-003}
 Let $f: B^n \to \IC^n$ be a $(K,K')$-quasiregular pluriharmonic mapping, $n>1$, such that $f(0)=0, \widetilde{\lambda}_f(0)=1$ and $\widetilde{\Lambda}_f(z)\leq \widetilde{\Lambda}$ for $z \in B^n$. Then $f$ is univalent on the ball $B^n(0,\rho)$ and $f(B^n(0,\rho))$ contains the ball $B^n(0, R)$, where $$\rho=\frac{\pi}{4(K+K'+\widetilde{\Lambda})} \ \  \text{and}\ \ R=\frac{\pi}{8(K+K'+\widetilde{\Lambda})}.$$
\end{thm}
\begin{pf}
 Let $z_1, z_2 \in B^n(0, \rho)$ be two fixed distinct points and
$z_1-z_2=|z_1-z_2| \theta$ for some $\theta\in \partial B^n$. Define the plurihamomic mapping
$$
\phi_\theta(z)=\left(f_z(z)-f_z(0)\right) \theta+\left(f_{\bar{z}}(z)-f_{\bar{z}}(0)\right) \bar{\theta}.
$$
Then, the definition of $\widetilde{\Lambda}_f(z)$ gives that
$$
|\phi_\theta(z)| \leq \widetilde{\Lambda}_f(z)+\widetilde{\Lambda}_f(0) \leq \widetilde{\Lambda}+K \widetilde{\lambda}_f(0)^{1 / n}+K^{\prime}=\widetilde{\Lambda}+K+K'
\ \ \text { for } z \in B^n. 
$$
Note that $\phi_\theta(0)=0$. By Lemma \ref{lemma3}, we obtain
$$
|\phi_\theta(z)| \leq \frac{4}{\pi}\left(\widetilde{\Lambda}+K+K' \right)|z| \quad \text { for } z \in B^n.
$$
We have
$$
\begin{aligned}
|f(z_2)-f(z_1)| =&\left|\int_{[z_1,z_2]}  f_z(z)dz+f_{\bar{z}}(z) d\bar{z} \right|\\ \geq & \left|\int_{[z_1, z_2]} f_z(0) d z+f_{\bar{z}}(0) d \bar{z}\right|
-\left|\int_{[z_1,z_2]}\left(f_z(z)-f_z(0)\right) d z+\left(f_{\bar{z}}(z)-f_{\bar{z}}(0)\right) d \bar{z}\right| \\
\geq & |z_2-z_1| \widetilde{\lambda}_f(0)-\int_{[z_1,z_2]}|\phi_\theta(z)| d s \\
> & |z_2-z_1|- \frac{4\left(\widetilde{\Lambda}+K+K'\right) \rho}{\pi}|z_2-z_1| =0.
\end{aligned}
$$
Thus $f(z_1)\neq f(z_2)$. This shows that $f$ is univalent in $B^n(0,\rho)$.

\hspace{-5mm}Now, let $z'\in \partial B^n(0,\rho)$. As $f(0)=0,$ we have
$$
\begin{aligned}
|f(z')|= & \left|\int_{[0,z']}f_z(z)dz+f_{\bar{z}}(z)d\bar{z}\right|\\
\geq & \left|\int_{[0, z']} f_z(0) d z+f_{\bar{z}}(0) d \bar{z}\right|
-\left|\int_{[0,z']}\left(f_z(z)-f_z(0)\right) d z+\left(f_{\bar{z}}(z)-f_{\bar{z}}(0)\right) d \bar{z}\right|\\
\geq & \widetilde{\lambda}_f(0) \rho- \int_{0}^\rho \frac{4(K+K'+\widetilde{\Lambda})r}{\pi} dr\\
= & \rho - \frac{2(K+K'+\widetilde{\Lambda})\rho^2}{\pi}=\frac{\pi}{8(K+K'+\widetilde{\Lambda})}=R.
\end{aligned}
$$
This shows that $f(B^n(0,\rho))$ contains the ball $B^n(0,R)$. This completes the proof of this theorem.
\end{pf}

Next, we establish a Bloch-type theorem for $(K,K')$-quasuregular pluriharmonic mappings.
\begin{thm}\label{Rohi-Vasu-P2-Theorem-004
}
Let $f: B^n \to \IC^n$ be a $(K,K')$-quasiregular pluriharmonic mapping, $n>1$, such that $\widetilde{\lambda}_f(0)=1$. Then $f(B^n)$ contains a schlicht ball of radius $b_f$, with
  $$\beta_f\geq \frac{\pi}{4(3K+2K')}.$$
\end{thm}
\begin{pf}
Let $f:B^n\to \IC^n$ be a $(K,K')$-quasiregular pluriharmonic mapping of the unit ball $B^n$ into $\IC^n$ with $\widetilde{\lambda}_f(0)=1$. Without loss of generality, we assume that $f$ is pluriharmonic on $\overline{B}^n$. We introduce the following functions 
$$M_1(r)=\max_{|z|\leq r}\widetilde{\lambda}_f(z)^{1/n},$$
and
$$\phi_1(r)=(1-r)M_1(r)$$
for $0\leq r\leq1.$
It is easy to see that $\phi_1(0)=M_1(0)=\widetilde{\lambda}_f(0)^{1/n}=1$ and $\phi_1(1)=0.$ Then there exist $r_0$ such that $\phi_1(r_0)=1$ and $\phi_1(r)<1$ for $r_0<r\leq1$. 
\vspace{2mm}

Also, since the set $\{z:|z|\leq r_0\}$ is compact, there exist $z_0$ such that $|z_0|\leq r_0$ and $ M_1(r_0)=\widetilde{\lambda}_f(z_0)^{1 / n}$, which implies
$$
\phi_1(r_0)=(1-r_0) M_1(r_0)=(1-r_0) \widetilde{\lambda}_f(z_0)^{1 / n}.
$$
Therefore,
\begin{equation}\label{Rohi-Vasu-P2-equation-011}
\left(1-r_0\right) \widetilde{\lambda}_f(z_0)^{1 / n}=1.
\end{equation}
Let $z\in B^n$ with $|z|=r \geq r_0$, then
\begin{equation}\label{Rohi-Vasu-P2-equation-011a}
(1-|z|)\widetilde{\lambda}_f(z)^{1/n}\leq (1-r)M_1(r)\leq 1.
\end{equation}  
In particular, we have
\begin{equation}\label{Rohi-Vasu-P2-equation-012}
\widetilde{\lambda}_f(z)\leq \widetilde{\lambda}_f(z_0) \quad\text{for} \ |z|=r_0.
\end{equation}
We consider following two cases,

\textbf{Case 1.} $r_0>0$: Fix a point $w_0$ with $0<|w_0|\leq r_0$ and assume that $\widetilde{\Lambda}_f(w_0)=|f_z(w_0)\theta+f_{\bar{z}}(w_0)\bar{\theta}|$ with $\theta\in \partial B^n.$ Define the function $\varphi$ by $$\varphi(\zeta)=f_z(\zeta w_0/|w_0|)\theta + f_{\bar{z}}(\zeta w_0/|w_0|)\bar{\theta} \quad \text{for}\ \zeta\in \ID.$$
Since $\varphi$ is harmonic, by the maximum modulus principle, there exists a point $\zeta'$ with $|\zeta'|=r_0$, such that 
   $$\widetilde{\Lambda}_f(w_0)=|\varphi(|w_0|)|\leq |f_z({\zeta' w_0}/{|w_0|})\theta+f_{\bar{z}}({\zeta' w_0}/{|w_0|})\bar{\theta}|. $$
Let $z_1=\zeta'w_0/|w_0|.$ Since $|z_1|=r_0$, by the definition of $(K,K')$-quasiregular pluriharmonic mappings and (\ref{Rohi-Vasu-P2-equation-012}), we have
 $$
 \begin{aligned}
   \widetilde{\Lambda}_f(w_0)\leq & \left|f_z(z_1)\theta+f_{\bar{z}}(z_1)\bar{\theta}\right|\\
    \leq & \widetilde{\Lambda}_f(z_1)\leq K\widetilde{\lambda}_f(z_1)^{1/n}+K'\leq  K\widetilde{\lambda}_f(z_0)^{1/n}+K'.
 \end{aligned} 
 $$ 
On the other hand, by the definition of $(K,K')$-quasiregular pluriharmonic mappings and (\ref{Rohi-Vasu-P2-equation-011}) with $\widetilde{\lambda}_f(0)=1,$ we have
  $$\widetilde{\Lambda}_f(0)\leq K\widetilde{\lambda}_f(0)^{1/n}+K'=K+K'=K(1-r_0)\widetilde{\lambda}_f(z_0)^{1/n}+K'.$$
This shows that 
\begin{equation}\label{Rohi-Vasu-P2-equation-013}
   \widetilde{\Lambda}_f(z)\leq  K\widetilde{\lambda}_f(z_0)^{1/n}+K' \quad \text{for} \ |z|\leq r_0.
\end{equation}
For $\xi \in B^n$, define
\begin{equation}\label{Rohi-Vasu-P2-equation-013a}
g(\xi)=z_0+\frac{(1-r_0)^n}{2}\xi \quad \text{and}\quad F(\xi)=2\left(f(g(\xi))-f(z_0)\right).
\end{equation} 
Then, it is easy to see that $$F(0)=0 \quad \text{and}\ \widetilde{\lambda}_F(0)=(1-r_0)^n \widetilde{\lambda}_f(z_0)=1.$$
If $|g(\xi)|\leq r_0$, from (\ref{Rohi-Vasu-P2-equation-013}) and (\ref{Rohi-Vasu-P2-equation-011}) we have
$$
\begin{aligned}
\widetilde{\Lambda}_F(\xi) = & (1-r_0)^n \widetilde{\Lambda}_F(g(\xi)) \\
\leq &  K(1-r_0)^n\widetilde{\lambda}_f(z_0)^{1/n}+(1-r_0)^nK'\\
\leq & K(1-r_0)\widetilde{\lambda}_f(z_0)^{1/n}+K'=K+K',
\end{aligned}
$$
and if $|g(\xi)|\geq r_0$, from (\ref{Rohi-Vasu-P2-equation-011a}), we obtain
\begin{align*}{\label{Rohi-Vasu-P2-equation-014}}
\widetilde{\Lambda}_F(\xi)= & (1-r_0)^n \widetilde{\Lambda}_F(g(\xi)) \\
\leq &  K(1-r_0)^n\widetilde{\lambda}_f(g(\xi))^{1/n}+(1-r_0)^nK'\\
\leq & K(1-r_0)\widetilde{\lambda}_f((g(\xi))^{1/n}+K'\\
= & K\left(\frac{1-r_0}{1-|g(\xi)|}\right)(1-|g(\xi)|)\widetilde{\lambda}_f((g(\xi))^{1/n}+K'\\ \vspace{2mm}
\leq & K\left(\frac{1-r_0}{1-|g(\xi)|}\right)+K'\leq \frac{K(1-r_0)}{1-r_0-(1-r_0)^n|\xi|/2}+K'\\ \vspace{2mm}
\leq & \frac{K(1-r_0)}{1-r_0-(1-r_0)|\xi|/2}+K'\\
= & \frac{2K}{2-|\xi|}+K'. \numberthis
\end{align*}

\textbf{Case 2.} $r_0=0$: Consider the functions $g$ and $F$ defined by (\ref{Rohi-Vasu-P2-equation-013a}) with $r_0=0$. Then $|g(\xi)|\geq r_0=0$ and it follows from (\ref{Rohi-Vasu-P2-equation-014}) that
 $$
 \widetilde{\Lambda}_F(\xi) = \frac{2K}{2-|\xi|}+K' \quad \text{for}\ \xi\in B^n.
 $$
 Therefore, we conclude that
 $$
\widetilde{\Lambda}_F(\xi) < 2K+K' \quad \text{for}\ \xi\in B^n.
 $$
 In particular, $\widetilde{\Lambda}_F(0)\leq K+K'.$
 \vspace{2mm}
 
 Now, applying Theorem \ref{Wang-2022} to the mapping $F$, we see that $F(B^n)$ contains a schlicht ball with the center $0$ and radius
       $$R'=\frac{\pi}{8(3K+2K')}.$$
Consequently, $f(B^n)$ contains a schlicht ball of radius 
$$R=\frac{\pi}{16(3K+2K')}.$$
This completes the proof.
\end{pf}

\noindent\textbf{Acknowledgement:}  The first named author thank SERB-CRG and the second named author thank CSIR for their support.


\begin{thebibliography}{99}
%\bibitem{Ahlfors-1938}  {\sc L.V. Ahlfors}, An extension of Schwarz’s lemma, {\it  Trans. Am. Math. Soc.} {\bf 43} (1938), 359--364. 
\bibitem{Ahlfors-Grunsky-1937}  {\sc L.V. Ahlfors} and {\sc H. Grunsky}, Uber die Blochsche Konstante, {\it  Math. Z. } {\bf 43} (1937), 671–-673.

\bibitem{Allu-Kumar-2024}{\sc V. Allu} and {\sc R. Kumar}, Landau-Bloch type theorem for elliptic and quasiregular harmonic mappings, {\it  J. Math. Anal. Appl. } (2024), https://doi.org/10.1016/j.jmaa.2024.128215
\bibitem{Bloch-1925}
{\sc A. Bloch}, Les théorèmes de M. Valiron sur les fonctions entières et la théorie de l'uniformisation, {\it  Ann.
Fac. Sci. Univ. Toulouse } {\bf 17} (1925), 1--22.

\bibitem{Bochner-1946}
{\sc S. Bochner}, Bloch's theorem for real variables, {\it Bull. Amer. Math. Soc.} {\bf 52 } (1946), 715--719.
 
 \bibitem{Bochner-Martin-1948}{\sc S. Bochner} and {\sc W. T. Martin}, Several complex variables, {\it  Princeton Math. Ser.} {\bf 10} (1948), p. 59.
 
 \bibitem{Bonk-1990}  {\sc M. Bonk}, On Bloch’s constant, {\it   Proc. Amer. Math. Soc.} {\bf 110} (1990), 889--894.
 
  \bibitem{Bconstant} {\sc H.H. Chen } and {\sc P.M. Gauthier},  On Bloch’s constant, {\it  J. Anal. Math.} {\bf 69} (1996), 275--291.
  
  \bibitem{Chen-Gauthier-2000} {\sc H.H. Chen} and {\sc P.M. Gauthier}, Bloch constants in several variables, {\it Trans. Amer. Math. Soc.} {\bf 353} (2000), 1371--1386. 
  
  \bibitem{Chen-Gauthier-2011} {\sc H.H. Chen} and {\sc P. M. Gauthier}, The Landau and Bloch theorem for planer harmonic and pluriharmonic mappings, {\it Proc. Amer. Math. Soc.} {\bf 139}(2) (2011), 583--595. 
  
    \bibitem{Chen-Li-Sahoo-Wang-2017} {\sc J. Chen}, {\sc P. Li}, {\sc S.K. Sahoo} and {\sc X. Wang}, On the Lipschitz continuity of certain quasiregular mappings between smooth Jordan domains, {\it Israel J. Math.} {\bf 220} (2017), 453--478.
  
  \bibitem{Chen-Ponnusamy-2020} {\sc S. Chen } and {\sc S. Ponnusamy}, On certain quasiconformal and elliptic mappings, {\it  J. Math. Anal. Appl.} {\bf 486} (2020), 1--16.
  
  \bibitem{Chen-Ponnusamy-Wang-2012-a} {\sc S. Chen}, {\sc S. Ponnusamy} and {\sc X. Wang}, Equivalent moduli of continuity, Bloch’s theorem for pluriharmonic mappings in $B^n$, {\it Proc. Indian Acad. Sci. Math. Sci.} {\bf 122} (2012), 583--595. 

 
  \bibitem{Chen-Ponnusamy-Wang-2012-b} {\sc S. Chen}, {\sc S. Ponnusamy} and {\sc X. Wang}, Landau–Bloch Constants for Functions in $\alpha$-Bloch Spaces and Hardy Spaces, {\it  Complex Anal. Oper. Theory} {\bf 6} (2012), 1025--1036. 

 
  \bibitem{Chen-Ponnusamy-Wang-2012-c} {\sc S. Chen}, {\sc S. Ponnusamy} and {\sc X. Wang}, Weighted Lipschitz continuity, Schwarz–Pick’s lemma and Landau–Bloch’s theorem for hyperbolic-harmonic mappings in $\IC^n$, {\it  Math. Model. Anal.} {\bf 18}(1) (2012), 66--79. 
   
 \bibitem{Chen-Ponnusamy-Wang-2015} {\sc S. Chen}, {\sc S. Ponnusamy} and {\sc X. Wang}, Stable geometric properties of pluriharmonic and biholomorphic mappings, and Landau–Bloch’s theorem, {\it Monatsh. Math.} {\bf 177} (2015), 33--51.   
 
 
  \bibitem{Chen-Ponnusamy-Wang-2021} {\sc S. Chen}, {\sc S. Ponnusamy} and {\sc X. Wang}, Remarks on Norm Estimates of the Partial Derivatives for Harmonic Mappings and Harmonic Quasiregular Mappings, {\it  J. Goem. Anal.} {\bf 31} (2021), 11051--11060. 
  
  %\bibitem{Clunie-Sheil-1984} {\sc J.G. Clunie} and {\sc T. Sheil-Small}, Harmonic univalent functions, {\it  Ann. Acad. Sci. Fenn. Ser. A. I.} {\bf 9} (1984), 3--25. 
  
   \bibitem{Eremenko-2000}{\sc A. Eremenko}, Bloch radius, normal families and quasiregular mappings, {\it Proc. Amer. Math. Soc.}, {\bf 128}(2), (2000), 557--560.
  
   \bibitem{Finn-Serrin-1958} {\sc R. Finn} and {\sc J. Serrin}, On the H\"{o}lder continuity of quasiconformal and elliptic mappings, {\it  Trans. Amer. Math. Soc.} {\bf 89} (1958), 1--15.

 
  \bibitem{Hahn-1973} {\sc K. T. Hahn}, Higher dimensional generalizations of the Bloch constant and their lower bounds, {\it Trans. Amer. Math. Soc.} {\bf 179} (1973), 263--274. 
  
  
  \bibitem{Harris-1977} {\sc L. A. Harris}, On the size of balls covered by analytic transformations, {\it Monatsh. Math.} {\bf 83} (1977), 9--23. 
  % \bibitem{Kalaj-2008} {\sc D. Kalaj}, Quasiconformal and harmonic mappings between Jordan domains, {\it  Math. Z.} {\bf 260} (2008), 237--252.
   
 % \bibitem{Kalaj-Ponnusamy-Matti-2014}{\sc D. Kalaj,  S. Ponnusamy} and {\sc M. Vuorinen}, Radius of close-to-convexity and full starlikeness of harmonic mappings, {\it Complex Var. Elliptic Equ.} {\bf 59} (2014), 539--552.

 %   \bibitem{Landau-1926} {\sc E. Landau}, Der Picard-Schottkysche Satz und die Blochsche Konstanten, {\it  Sitzungsber.Preuss. Akad. Wiss. Berlin Phys.-Math. Kl.} (1926), 467--474.

 \bibitem{Landau1} {\sc E. Landau},  Uber die Bloehsehe Konstante und zwei verwandte Weltkonstanten, {\it Math. Z.} {\bf 30} (1929), 608--634.

  \bibitem{Liu-Ponnusamy-2022} {\sc M.S. Liu} and {\sc S. Ponnusamy},  Bloch and Landau type theorems for pluriharmonic mappings, {\it  Internat. J. Math.} {\bf 33}(7) (2022), 2250053.
  
  \bibitem{Marden-Rickman-1974} {\sc A. Marden} and {\sc S. Rickman}, Holomorphic mappings of bounded distortion, {\it  Proc. Amer. Math. Soc.} {\bf 46} (1974), 225--228.

    \bibitem{Nirenberg-1953} {\sc L. Nirenberg}, On nonlinear elliptic partial differential equations and H\"{o}lder continuity, {\it  Commun. Pure. Appl. Math.} {\bf 6} (1953), 103--156.

   \bibitem{Poletsky-1985} {\sc E. A. Poletsky}, Holomorphic quasiregular mappings, {\it Proc. Amer. Math. Soc.} {\bf 92} (1985), 235--241.

  \bibitem{Rudin-1980}{\sc W. Rudin}, Function theory in the unit ball of $\IC^n$, { Springer-Verlag, New York, Heidelberg, Berlin}, 1980.

\bibitem{Sakaguchi-1956} {\sc K. Sakaguchi}, On Bloch's theorem for several complex variables, {\it  Sci. Rep. Tokyo Kyoiku Daigaku Sect.} {\bf 5} (1956), 149--154.

\bibitem{Takahashi-1951} {\sc S. Takahashi}, Univalent mappings in several complex varibles, {\it  Ann. of Math.} {\bf 53} (1951), 464--471.

  \bibitem{Wang-Yang-Liu-2022} {\sc X. Wang}, {\sc Y. Yang} and {\sc M. S. Liu }, The Landau-Bloch type theorems for $K$-quasiregular pluriharmonic mappings, {\it Monatsh. Math.} {\bf 198}(1) (2022), 189--209. 
 
\bibitem{Wu-1967} {\sc H. Wu}, Normal families of holomorphic mappings, {\it  Acta Math.} {\bf 119} (1967), 193--233.

\bibitem{Xu-Liu-2020} {\sc Z.F. Xu} and {\sc M.S. Liu}, On pluriharmonic $\nu$-Bloch-type mappings and hyperbolic-harmonic mappings, {\it  Monatsh. Math.} {\bf 192} (2020), 965--978. 

\bibitem{Yu-Gu-1997} {\sc Yu Yi-Sheng} and {\sc Gu Dun-he}, A note on a lower bound for the smallest singular value, {\it  Linear Algebra Appl.} {\bf 253} (1997), 25--38.

 
\end{thebibliography}
\end{document}